\documentclass[12pt, leqno]{article}

\usepackage[english]{babel}

\usepackage{amsthm, amssymb, amsmath}

\theoremstyle{plain}
 \newtheorem{thm}{Theorem}
 
 \newtheorem{lem}[thm]{Lemma}
 \newtheorem{prop}[thm]{Proposition}
 
 \theoremstyle{definition}
 \newtheorem{defn}[thm]{Definition}
 \theoremstyle{remark}
 \newtheorem{rem}[thm]{Remark}

\begin{document}

\title{On Algebraic Shift Equivalence of Matrices  over Polynomial Rings}

\author{by \large SHENG  CHEN }

\maketitle

{\small\it (Department of Mathematics,   Harbin Institute of
Technology,
Harbin 150001,   P.R.China)}\\

\footnotetext{\hspace{-0.5cm} $^*$\hspace*{1mm} 
 {\it E-mail address:} {schen@hit.edu.cn.} }

\begin{center}
{\bf Abstract\\[1ex]}
\end{center}
The paper studies  algebraic  shift equivalence of matrices over
$n$-variable polynomial rings over a principal ideal domain
$D$($n\leq 2$). It is proved that in the case $n=1$, every
non-nilpotent matrix over $D[x]$ is
 algebraically strong shift equivalent to a nonsingular matrix. In the case
 $n=2$,
 an example of non-nilpotent matrix  over $\mathbb{R}[x,y,z]=\mathbb{R}[x][y,z]$,
 which can not be algebraically shift equivalent to a nonsingular
matrix, is
 given.\\

 \noindent 2000 AMS Classification:  Primary   15A54,   Secondary 15A23, 13C10, 37B10\\
 Key words: \hspace{1mm}
 matrix, polynomial ring, full rank factorization, algebraic  shift
 equivalence, projective module

\section{Introduction}
Let $A$ and $B$ be two square matrices  of possibly  different order
over a semi-ring  $ R$ with 1. We say $A$ and $B$ are
\emph{elementarily
 strong shift equivalent  over  $ R$}  if there exist  matrices $U,V$ over
$ R$ such that
$$A=UV , VU=B$$
In this case we write $(U,V): A\equiv B$ or simply $ A\equiv B$.

 We say $A$ and $B$ are
\emph{algebraically strong shift equivalent over  $ R$ of lag $l$}
and we write $A\approx B$(lag $l$)  if there exists a sequence of
$l$ elementary equivalence from $A$ to $B$:
$$(U_1,V_1): A=A_0\equiv A_1,(U_2,V_2): A_1\equiv A_2,\cdots,(U_l,V_l): A_{l-1}\equiv
A_l=B$$ Say that $A$ is \emph{algebraically strong shift equivalent}
 to $B$ and write $A\approx B$ if $A\approx B$(lag $l$) for some   $l\in \mathbb{N}$.

In symbolic dynamics,  the case $R=\mathbb{Z}_+=\mathbb{Z}\cap
[0,\infty)$ is the most interesting (cf \cite{W73} or  \cite{Lind95}
). However, the problem of the decidability of algebraic strong
shift equivalence of matrices over $R=\mathbb{Z}_+$ is still open.

Let $A$ and $B$ be matrices over a semi-ring  $ R$ with 1. If there
exists a pair $U$ and $V$  of matrices over $R$ such that
$$AU=UB,VA=BV,A^l=UV,B^l=VU$$
then we say that $A$ and $B$ are \emph{ algebraically shift
equivalent of lag $l$ over $R$}. We denote this situation by
$(U,V):A\sim B $(lag $l$). We say $A$ is \emph{algebraically shift
equivalent}
 to $B$ (and write $A\sim B$) if $A\sim B$(lag $l$) for some  $l$.

It is easy to check that  \emph{algebraically strong shift
equivalence} implies  \emph{algebraically shift equivalence}. As is
well-known,   K.H. Kim and F.W. Roush showed  in \cite{Kim3} that
the problem of algebraic shift equivalence of matrices over
$R=\mathbb{Z}_+$ is decidable.

We say that a domain $D$  has \emph{property  NSSEN} if \emph{every
non-nilpotent matrix $A$  over  $D$ is algebraically strong shift
equivalent  to some nonsingular matrix.}

Similarly, we say that a domain $D$  has \emph{property  NSEN} if
\emph{every non-nilpotent matrix $A$ over $D$ is algebraically shift
equivalent to some nonsingular matrix.}

In \cite{EF81},  E.G. Effros  showed that a principal ideal domain
$R$ has  \emph{property  NSSEN}.  In \cite{trans},  M. Boyle and D.
Handelman showed  that  a commutative domain $R$ admitting a
non-free finitely generated projective module does \emph{not}  have
\emph{property NSEN}. It is well-known that the Quillen-Susslin
Theorem(cf. \cite{q76} or \cite{s76}) says  that  a projective
module over  the ring of $n$-variable polynomials  over a principal
ideal domain $D$  is always free. So it is natural to ask the
following problem.

\emph{Does the ring of $n$-variable  ($n\in \mathbb{N}$) polynomials
over a principal ideal domain $D$  has property NSSEN or property
NSEN}?

 This paper aims  to answer the problem positively
  in the case $n=1$ through proving the existence of full
rank factorization for any matrix over the ring of univariate
polynomials over $D$, and answer the question negatively  in the
case $n=2$ by a counterexample.

The main   results of the paper are as follows.\\

{\thm \label{th1}
 Let $D$ be a principal ideal domain and $A$ be a square  matrix  of order $n$ over
 $D[x]$.
If $A$ is not nilpotent and
$l=min\{k\in\mathbb{N}|rank(A^k)=rank(A^{k+1})\}$, then there is an
algebraic strong shift equivalence of lag $l$ from  $A$ to a
nonsingular matrix. In other words, $D[x]$ has property NSSEN. }


{\thm \label{th2} Let
$$A= \left(
    \begin{array}{ccc}
      0 &z  & -y \\
      -z &0 & x\\
      y & -x & 0 \\
    \end{array}
  \right)
$$
Then $A$  can not be algebraically shift equivalent to a nonsingular
matrix  over $\mathbb{R}[x,y,z]=\mathbb{R}[x][y,z]$. In other words,
$\mathbb{R}[x,y,z]$ does not
has property NSEN. }\\

 The organization of the
paper is as follows. In section 2 we prove the existence of full
rank factorization of matrices over univariate polynomials. In
section 3 we give the proof of the main results.

\section{Existence of Full Rank Factorization}

The aim of this section is to prove the following result, which will
be used in section 3.  The result is a slight generalization of the
corresponding  result in page 109 of \cite{79} and our proof is
based on Theorem 2 in \cite{82} and some results in \cite{79}.

 \begin{prop}\label{pr2}  Let $D$ be a principal ideal domain and $A$ be a matrix of size $m\times n$ with    rank   $r$  over
 $D[x]$.
Then there exists a full rank factorization of $A$ over
 $D[x]$:
$$A=PQ$$
where $P$ and $Q$ are matrices  of size $m\times r$  and $r\times n$
over $D[x]$ respectively.
\end{prop}

To prove the proposition above, we  first recall some facts about
matrices over domains(cf. \cite{BR84}), which are assumed to be
commutative throughout the paper, and prove some lemmas.

Let $D$ be a principal ideal domain. Throughout this section the
matrices will be over $D[x]$, which is a unique factorization domain
(UFD). Similar to the case of fields, every non-zero matrix has a
rank defined by its maximal order of non-zero  minors. A square
matrix $A$ is invertible if and only if its determinant is a unit in
$D$. Let $L$ and $U$ be $m\times r$ and $r \times n$ matrices over
$D[x]$ respectively. Then the Binet-Cauchy theorem says that
$\wedge^t(L)\wedge^t(U)=\wedge^t(LU)$, where   $\wedge^t(L)$ and
$\wedge^t(U)$ are the $t$-th compound matrix of $L$ and $U$
respectively for $t\leq min \{m,r,n\}$.

 {\defn
Let $C$ be an $m\times n$ matrix over $D[x]$ with rank $m$. We say
that  $C$ is  \emph{minor left prime (MLP)} if 1 is the greatest
common polynomial  divisor (gcd) of  all the $m\times m$ minors. Say
that  $C$ is \emph{factor left prime (FLP)} if in any polynomial
matrix factorization $C=C_1C_2$, where $C_1$ is square,  $C_1$
should  be an invertible matrix.


{\lem\label{lem1} (cf. Theorem 2 in \cite{82})
 Let $D$ be a principal ideal domain and  $A$ be an $m\times n$  matrix ($m\leq n$) of rank $m$ with entries in
 $D[x]$.
 Let $d(x)$ be  the  greatest common  divisor of all $m$-th order
 minors of $A$. Then $A$ has a factorization as $A=LU$ with det
 $L=d(x)$. If $D$ is an Euclidean domain, then we have algorithm to find the factorization. }\\

 {\rem In \cite{82},  the  lemma  above is proved  only  for  Euclidean domains. However,  the authors of \cite{82} remarked  in page
 656 that the result is also true for  principal ideal domains.}

 {\lem\label{lem2} For a matrix $A$ $m\times n$  matrix ($m\leq n$) of rank $m$ over $D[x]$, $A$ is MLP
if and only if $A$ is FLP.}\\

\begin{proof}  Suppose that $A$ is MLP and  $A=LU$ is a factorization.  If $\det(L)=f(x)$, then $\bigwedge ^m(A)=\bigwedge^m
(L)\bigwedge^m (U)$ and thus $f(x)|gcd(\bigwedge ^m(A))=1$.  So $L$
should be invertible. We have MLP $\Rightarrow$ FLP.  Next,  it
follows from  Lemma \ref{lem1} that FLP $\Rightarrow$ MLP. The proof
is complete.
\end{proof}

{\rem Lemma \ref{lem2} is a slight generalization of Theorem  3 in
\cite{79}, which says that  MLP and FLP are equivalent for matrices
over the ring of bivariate polynomials over fields.}

{\lem \label{lem3} A matrix  $C$ of size $m \times n$  with rank $m$
over $D[x]$ is MLP if and only if there exist $d_j\in D$ and
matrices $Z_j$ such that
$$CZ_j=d_jI_m$$ where $0 \leq j \leq s$, $s \in \mathbb{N}\cup \{0\}$ and $gcd(d_0,d_1,d_2,\cdots,d_s)=1$.}\\

\begin{proof}
Denote    by $C_{i_1,i_2,\cdots,i_m}$ the  $m\times m$ submatrix
from $i_1,i_2,\cdots, i_m$ columns  of $C$, by
$C_{i_1,i_2,\cdots,i_m}^*$ the  adjoint matrix of
$C_{i_1,i_2,\cdots,i_m}$,  and by $\Delta_{i_1i_2\cdots i_m}$ the
determinant of $C_{i_1,i_2,\cdots,i_m}$. Then
$$C_{i_1,i_2,\cdots,i_m}C_{i_1,i_2,\cdots,i_m}^*=\Delta_{i_1i_2\cdots i_m}I_m$$
Define an $l$ by $m$ matrix  $Z_{i_1i_2\cdots i_m}$ as follows: its
$i_1,i_2,\cdots, i_m$-th rows are  the $1, 2,\cdots, m$-th rows of
$C_{i_1,i_2,\cdots,i_m}^*$  and all other rows are 0. Then we have
$$CZ_{i_1i_2\cdots i_m}=\Delta_{i_1i_2\cdots i_m}I_m$$

Suppose that $C$ be MLP. Let $F$ be the quotient field of $D$. Then
$D[x]$ is a subring of $F[x]$, which is a principal ideal domain.
Since 1 is the gcd of $\{ \triangle_{i_1i_2\cdots i_m}| 1\leq
i_1\leq i_2\leq\cdots \leq i_m\leq n\}$, there exist
$a_{0,i_1i_2\cdots i_m} \in D[x]$ such that
$$ \sum_{1\leq
i_1\leq i_2\leq\cdots \leq i_m\leq n } a_{0,i_1i_2\cdots
i_m}\triangle_{i_1i_2\cdots i_m}= d_0 \in D \setminus\{0\}$$

If $d_0=1$, then let $s=0$. Otherwise, suppose that
$$d_0=\prod_{j=1}^{s} p_j^{t_j}$$ where $p_j$ are distinct prime factors of $d_0$ in
$D$, $t_j\in \mathbb{N}$,  and $s\in \mathbb{N}$. Since the gcd of
$\{ \triangle_{i_1i_2\cdots i_m}| 1\leq i_1\leq i_2\leq\cdots \leq
i_m\leq n\}$ is 1, for any $j\in \{1,2, \cdots s\}$,  there exists
at least one of them,
 say  $\triangle_{j_1j_2\cdots j_m}$, is
not divisible by $p_j$. Let $d_0+ \triangle_{j_1i_2\cdots j_m}
=d_j$.
 Define

$$ a_{j,i_1i_2\cdots i_m}=\{ \begin{array}{lcc}
                                 a_{0,i_1i_2\cdots i_m} &
 \text{ if\, \,} &  i_1i_2\cdots i_m\neq {j_1j_2\cdots j_m} \\
                                  a_{0,i_1i_2\cdots i_m} +1& \text{ if\, \,} & i_1i_2\cdots i_m=
{j_1j_2\cdots j_m}
                                \end{array}
$$
Let $$Z_j=\sum_{1\leq i_1\leq i_2\leq\cdots \leq i_m\leq n}
a_{j,i_1i_2\cdots i_m}Z_{i_1i_2\cdots i_m}
$$  Then
$$CZ_j=   \sum_{1\leq
i_1\leq i_2\leq\cdots \leq i_m\leq n} a_{j,i_1i_2\cdots
i_m}\triangle_{i_1i_2\cdots i_m}I_m=( d_0+ \triangle_{j_1j_2\cdots
j_m}) I_m=d_jI_m$$
 Obviously $gcd(d_0,d_1,d_2,\cdots,d_s)=1$.
The proof is complete. \end{proof}

{\rem The proof  of Lemma \ref{lem3} is adapted from  that  of
Theorem 2 in \cite{79}.  }

 {\lem \label{lem4} Suppose that $A_{11},A_{12}$ and
$A_{21}$ are three matrices over $D[x]$ such that
$A_{21}A_{11}^{-1}A_{12}$ is a matrix over $D[x]$. If
$C=(A_{11},A_{12})$ is MLP, then $A_{21}A_{11}^{-1}$ is
a matrix over  $D[x]$.}\\

\begin{proof} By Lemma \ref{lem3}, we can prove the result using almost the same method in the proof of
the corollary to Theorem 2 in \cite{79}. The detail is omitted.
\end{proof}


 Now we can prove Proposition \ref{pr2}.
\begin{proof} It suffices to consider the case $A$ being an $m\times n$ polynomial matrix
with rank $r<min(m,q)$. By an appropriate interchange of rows and
columns we can always  assume that $$ A=\left(\begin{array}{cc}
     A_{11} & A_{12}\\
    A_{21} & A_{22}
   \end{array}
\right) $$ where  $A_{11}$ is a nonsingular matrix of order $r$. By
Lemma \ref{lem1} and Lemma  \ref{lem2}, we can suppose that
$(A_{11},A_{12})=L(U_{11},U_{12})$, where $ (U_{11},U_{12})$ is MLP.
Then
$$A=\left(\begin{array}{cc}
     L & 0\\
   0 & I_{m-r}
   \end{array}
\right)\left(\begin{array}{cc}
     U_{11} & U_{12}\\
    A_{21} & A_{22}
   \end{array}
\right) $$ Note that
$$A_{22}=A_{21}U_{11}^{-1}U_{12}$$
By Lemma \ref{lem4}, $A_{21}U_{11}^{-1}$ is a matrix over $D[x]$.
Obviously, we have
$$\left(\begin{array}{cc}
     U_{11} & U_{12}\\
    A_{21} & A_{22}
   \end{array}
\right)=A_1A_2$$where
$$A_1=\left(\begin{array}{c}
   I_r\\
    A_{21}U_{11}^{-1}
   \end{array}
\right), A_2=\left(U_{11}, U_{12}\right)$$ Let
$$A_3=\left(\begin{array}{c}
   L\\
    A_{21}U_{11}^{-1}
   \end{array}
\right)
$$
 Then
$$A=A_3A_2$$
is a full rank factorization.
 \end{proof}
{\rem Although the proof  above is almost the same as  that  in page
108-109 in \cite{79}, it is included here for the readers'
convenience. }

\section{Proof of The Main Results}

We first prove  Theorem \ref{th1}.
\begin{proof}
Let $A\in M_{n}(D[x])$. Suppose that
$$A=B_1C_1,C_1B_1=B_2C_2,C_2B_2=B_3C_3,\cdots$$
is  a sequence of full rank factorizations,  i.e.,  $B_iC_i$ are
full rank factorizations of
 $C_{i-1}B_{i-1}$, for $i=2,3,\cdots $.  If
$C_iB_i$ is $p\times p$ and has rank $q<p$, then the size of
$C_{i+1}B_{i+1}$ will be $q\times q$. That is, the size of
$C_{i+1}B_{i+1}$ must be strictly smaller than that of $C_iB_i$ when
$C_iB_i$ is singular.
  It follows that there eventually must be a pair of factors $B_k$
  and $C_k$, such that $C_kB_k$ is   either nonsingular or is
  zero. Let $l$ be  the first integer for which this occurs.
  Write
  $$A^l=(B_1C_1)^l=B_1(C_1B_1)^{l-1}C_1=\cdots=B_1B_2\cdots
  B_{l-1}(B_lC_l)
  C_{l-1}C_{l-2}\cdots C_1$$
 and $$A^{l+1}=B_1B_2\cdots
  B_{l-1} B_{l}(C_lB_l)C_l
  C_{l-1}C_{l-2}\cdots C_1$$

 If $C_lB_l=0$, then $A^{l+1}=0$.  So  $A$ is algebraically strong shift
equivalent to a nonsingular matrix $C_lB_l$ if $A$ is not nilpotent.
Assume that $C_lB_l$ is nonsingular. If $B_l$ is $p\times r$ and
$B_l$ is $r\times p$, then rank $B_lC_l=r$. Since $C_lB_l$ is
$r\times r$ and nonsingular, it follows that
$rank(C_lB_l)=r=rank(B_lC_l)$. Note that $B_i$ and $C_i$ are of full
column rank and full row rank respectively for $i=1,2,3,\cdots$. It
follows that
$$rank(A^{l+1})=rank(C_lB_l)=rank(B_lC_l)=rank(A^{l})$$

Finally, it is easy to check that
$l=min\{k\in\mathbb{N}|rank(A^k)=rank(A^{k+1})\}$.
 The proof is now complete.
\end{proof}

{\rem The proof  above is adapted from  that  of a theorem about the
existence  of Drazin inverse of matrices over fields in \cite{Wang}.
}

To prove Theorem \ref{th2}, we need the following two lemmas.

{\lem (cf. \cite{trans}) Let $R$ be a commutative domain admitting a
finitely generated projective module $P$ that is not free. From a
module $Q$ such that $P\oplus Q \simeq R^n$, let $e:  R^n
\rightarrow R^n $ be the endomorphism that projects onto $P$. Then
any matrix representation for $e$ is not algebraically shift
equivalent  to a nonsingular matrix over $R$.
  }

{\lem\label{weibel} (cf. Example 1.2.2 in \cite{weibel}) Let
$R=\mathbb{R}[x,y,z]/(x^2+y^2+z^2-1)$  and $\sigma$ be the following
homomorphism induced by  the unimodular  row $\alpha=(x,y,z)$:

$$ \sigma: R^3 \rightarrow R:
\begin{array}{ccc}
  \left(\begin{array}{c}
     a \\
     b \\
     c
   \end{array}
   \right)& \mapsto &ax+by+cz
\end{array}
$$
Then $Q\oplus R\simeq R^3$ and $Q=ker(\sigma)$ is a stably free
(projective)
module which is not free.}\\

 Now we can give the proof of Theorem \ref{th2}.
\begin{proof}
Let $ \varphi
:R=\mathbb{R}[x,y,z]\rightarrow\mathbb{R}[x,y,z]/(x^2+y^2+z^2-1)$ be
the natural  surjective ring homomorphism. Note that
$$A^3=-(x^2+y^2+z^2)A$$
  We have
   $$ \varphi{(A^4)}=  \varphi(-A^2)=\overline{\left(    \begin{array}{ccc}
                     1& 0& 0 \\
                    0 & 1 & 0 \\
                    0& 0 & 1
                    \end{array}
           \right)}-\overline{\left(    \begin{array}{c}
                      x  \\
                      y  \\
                      z
                    \end{array}
           \right)
           \left( \begin{array}{ccc}
                      x  &y&
                      z
                    \end{array}
           \right)}$$ and thus $\varphi{(A^4)}$ is an idempotent matrix over $\mathbb{R}[x,y,z]/(x^2+y^2+z^2-1)$.

 Note that $im \varphi(A^4) \simeq ker(\sigma)$, where $\sigma$ be the
module homomorphism induced by  the unimodular  row $\alpha=(x,y,z)$
 in Lemma \ref{weibel}. Suppose that $A$ is algebraically shift
equivalent to  a nonsingular $2\times 2$ matrix, say $B$. Then,  we
have an algebraic shift equivalence between  $\varphi(A^4)$ and
$\varphi(B^4)$ and thus
det$(I_2-t\varphi(B^4))$=det$(I_3-t\varphi(A^4))$ =$(1-t)^2$. So
$\varphi(B^4)$ is   a nonsingular $2\times 2$ matrix over
$\mathbb{R}[x,y,z]/(x^2+y^2+z^2-1)$.  A contradiction to Lemma
\ref{weibel}  arises. The proof is complete.
\end{proof}


\subsection*{Funding}
National Natural Science Foundation of China(10526016); National
Natural Science Foundation of China(10571033); Development Program
for Outstanding Young Teachers in Harbin Institute of Technology
(HIT.2006.54); Science Research Foundation in Harbin Institute of
Technology (HITC200701).

\subsection*{Acknowledgments}
The author would like to dedicate this paper to Prof. Hong YOU for
his introduction of him to the field of classical groups and
algebraic K-theory.

\bibliographystyle{amsplain}
\bibliography{xbib}
\end{document}